\documentstyle[12pt]{article}
\newcommand{\oid}[3]{\mbox{${\cal #1}_{#2}^{#3}$}}
\newcommand{\idn}[3]{\mbox{${\bf #1}_{#2}^{#3}$}}
\newcommand{\idrel}[3]{\mbox{${\cal #1} \stackrel{1}{#2} {\cal #3}$}}

\begin{document}
\newtheorem{theorem}{Theorem}[section]
\newtheorem{prop}{Proposition}[section]
\newtheorem{lemma}{Lemma}[section]
\newtheorem{cor}{Corollary}[section]

\title{Compositions of operator ideals and their regular hulls}
\author
{Frank Oertel\\
Zurich}
\maketitle

\begin{abstract}
{\noindent}Given two quasi-Banach ideals \oid{A}{}{} and \oid{B}{}{}
we investigate the regular hull of their composition -  
$(\oid{A}{}{} \circ \oid{B}{}{})^{reg}$. In concrete situations
this regular hull appears more often than the composition itself. As a first
example we obtain a description for the regular hull of the nuclear 
operators which is a "reflected" Grothendieck representation:\\
$\oid{N}{}{reg} \stackrel{1}{=} \oid{I}{}{} \circ \oid{W}{}{}$ (theorem 2.1).    
Further we recognize that the class of such ideals leads to interesting
relations concerning the question of the accessibility of (injective)
operator ideals.\\
 
{\noindent}{\it{Key words and phrases:}} accessibility, factorization, Banach spaces, 
operator ideals\\

{\noindent}{\it{1991 AMS Mathematics Subject Classification:}} primary 46M05, 47D50; 
secondary 47A80.  
\end{abstract}

\section{Introduction and notations}

In this paper we only deal with Banach spaces and most of our notations and definitions 
concerning Banach spaces and operator ideals are standard and can be
found in the detailed monographs \cite{df} and \cite{p1}. However, if 
$(\oid{A}{}{}, \idn{A}{}{})$ and $(\oid{B}{}{}, \idn{B}{}{})$ are given
quasi-Banach ideals,
we will use the shorter notation $(\oid{A}{}{d}, \idn{A}{}{d})$ for the dual
ideal (instead of $(\oid{A}{}{dual}, \idn{A}{}{dual})$) and the abbreviation 
\idrel{A}{=}{B} for the equality 
$(\oid{A}{}{}, \idn{A}{}{}) = (\oid{B}{}{}, \idn{B}{}{})$. The inclusion 
$(\oid{A}{}{}, \idn{A}{}{}) \subseteq (\oid{B}{}{}, \idn{B}{}{})$ is often
shortened 
by \idrel{A}{\subseteq}{B}, and if $T : E \longrightarrow F$
is an operator, we indicate that it is a metric injection by writing 
$T : E \stackrel{1}{\hookrightarrow} F$.\\
{\noindent}Each section of this paper includes the more special terminology which
is not so common.

{\noindent}Besides the maximal Banach ideal $(\cal L, \| \cdot \|)$
we will mainly be concerned with $(\oid{G}{}{}, \| \cdot \|)$ (approximable
operators), $(\oid{K}{}{}, \| \cdot \|)$ (compact operators),
$(\oid{W}{}{}, \| \cdot \|)$ (weakly compact operators), 
$(\oid{N}{}{}, \idn{N}{}{})$ (nuclear operators),
$(\oid{I}{}{}, \idn{I}{}{})$ (integral operators),
$(\oid{P}{p}{}, \idn{P}{p}{})$ (absolutely $p$-summing operators),
$1 \leq p \leq \infty, \frac{1}{p} + \frac{1}{q} = 1$
and $(\oid{L}{\infty}{}, \idn{L}{\infty}{}) \stackrel{1}{=} 
(\oid{P}{1}{\displaystyle\ast}, \idn{P}{1}{\displaystyle\ast})$.\\
{\noindent}Since it is important for us, we recall the notion 
of the conjugate operator ideal (cf. \cite{glr}, \cite{jo}): \\
let $(\oid{A}{}{}, \idn{A}{}{})$
be a quasi-Banach ideal. Let $\oid{A}{}{\Delta}(E, F)$ be the set of all $T \in {\cal L}(E, F)$
for which 
\begin{center}
$\idn{A}{}{\Delta}(T) : = \sup\{tr( TL) \mid L \in {\cal F}(F, E), \idn{A}{}{}(L) \leq 1 \} < \infty $.
\end{center}
Then a Banach ideal is obtained. It is called the {\it{conjugate ideal}} of
$(\oid{A}{}{}, \idn{A}{}{})$.\\

{\noindent}Denote for given Banach spaces $E$ and $F$
\begin{center}
FIN$(E) : = \{M \subseteq E\mid M\in $ FIN$\}$ \hspace{0.2cm}and\hspace{0.2cm}
COFIN$(E) : = \{L\subseteq E\mid E/L \in $ FIN$\}$,       
\end{center}
where FIN stands for the class of all finite dimensional Banach spaces.\\
A quasi-Banach ideal $(\oid{A}{}{}, \idn{A}{}{})$ is called {\it{right-accessible}}, if for all
$(M, F) \in$ FIN $\times$ BAN, operators $T \in {\cal L}(M, F)$ and $ \varepsilon > 0$ 
there are $N \in$ FIN$(F)$ and $S \in {\cal L}(M, N)$ such that $ T = J_N^F S$
and $\idn{A}{}{}(S) \leq (1 + \varepsilon) \idn{A}{}{}(T) $.
It is called {\it{left-accessible}}, if for all $(E, N) \in$ BAN $\times$ FIN,
operators $T \in {\cal L}(E, N)$ and $\varepsilon > 0$ there are $L \in$ COFIN$(E)$
and $S \in {\cal L}(E/L, N)$ such that $T = S Q_L^E $ and 
$\idn{A}{}{}(S) \leq (1 + \varepsilon) \idn{A}{}{}(T) $.
A left- and right-accessible ideal is called {\it{accessible}}. 
$(\oid{A}{}{}, \idn{A}{}{})$ is {\it{totally accessible}}, if for every finite rank operator
$T \in {\cal F}(E, F)$ between Banach spaces and $\varepsilon > 0$ there are
$(L, N) \in$ COFIN$(E) \times$ FIN$(F)$ and $S \in {\cal L}(E/L, N)$ such that 
$T = J_N^F S Q_L^E$ and $\idn{A}{}{}(S) \leq (1 + \varepsilon) \idn{A}{}{}(T) $. \\
Every injective quasi-Banach ideal is right-accessible (every surjective ideal is
left-accessible) and, if it is left-accessible, it is totally accessible.\\

\section{Compositions of operator ideals and applications to nuclear operators}

Let (\oid{A}{}{}, \idn{A}{}{}) be a $p$-Banach ideal and 
(\oid{B}{}{}, \idn{B}{}{}) be a $q$-Banach ideal $(0 < p, q \leq 1)$. Then:

\begin{lemma} 
$(\oid{A}{}{} \circ \oid{B}{})^{reg} \stackrel{1}{\subseteq} 
\oid{A}{}{reg} \circ \oid{B}{}{inj}$.
\end{lemma}

{\noindent}{\sc PROOF}: Let $E, F$ be Banach spaces, $\varepsilon > 0$ and $ T \in 
(\oid{A}{}{} \circ \oid{B}{}{})^{reg}(E, F)$. Then there are a Banach space $G$,
operators $R \in \oid{A}{}{}(G, F'')$ and $S \in \oid{B}{}{}(E, G)$ such that
$j_FT = RS$ and \\$\idn{A}{}{}(R) \idn{B}{}{}(S) < (1 + \varepsilon) 
(\idn{A}{}{} \circ \idn{B}{}{})^{reg}(T)$. Let $C$ be the (closed) range of 
$j_F : F \stackrel{1}{\hookrightarrow} F''$.  
Then $G_0 : = R^{-1}(C)$ is a closed subspace of $G$. Let $S_0 \in 
\oid{L}{}{}(E, G_0)$ be defined by $S_0x : = Sx$ $(x \in E)$. Then
$J_{G_0}^{G}S_0 = S \in \oid{B}{}{}(E, G)$. Hence $S_0 \in \oid{B}{}{inj}(E,G_0)$
and $\idn{B}{}{inj}(S_0) \leq \idn{B}{}{}(S)$. 
Now let $\gamma : C \longrightarrow F$ be defined canonically and 
$\gamma_0$ be the restriction of $\gamma$ to $C_0$ where $C_0$ is
the closure of $RJ_{G_0}^G(G_0)$. Setting $V : = \gamma_0R_0$ with 
$R_0 : G_0 \longrightarrow C_0$, $R_0z : = Rz$ $(z \in G_0)$
the construction implies that $j_FV = RJ_{G_0}^G \in \oid{A}{}{}(G_0, F'')$.
Hence $V \in \oid{A}{}{reg}(G_0, F)$, $\idn{A}{}{reg}(V) \leq 
\idn{A}{}{}(R)$ and $j_FT = RS = (RJ_{G_0}^G)S_0 = j_FVS_0$. It follows
that $T = VS_0 \in \oid{A}{}{reg} \circ \oid{B}{}{inj}(E, F)$ and 
$\idn{A}{}{reg}(V) \idn{B}{}{inj}(S_0) \leq 
\idn{A}{}{}(R) \idn{B}{}{}(S) < (1 + \varepsilon) (\idn{A}{}{} 
\circ \idn{B}{}{})^{reg}(T)$ and the proof is finished. 
${}_{\displaystyle\Box}$\\

\begin{cor}
Let $0 < p, q \leq 1, (\oid{A}{}{}, \idn{A}{}{})$ be a $p$-Banach ideal and 
$(\oid{B}{}{}, \idn{B}{}{})$ be a $q$-Banach ideal. 
Then $\oid{B}{}{reg} \stackrel{1}{\subseteq} \oid{B}{}{inj}$. If in addition $\oid{B}{}{}$ is
injective then $(\oid{A}{}{} \circ \oid{B}{}{})^{reg} \stackrel{1}{=}
\oid{A}{}{reg} \circ \oid{B}{}{}$.
\end{cor}
  
{\noindent}Next we show that lemma 2.1 (and the notion of accessibility)
yields a description of the regular hull of nuclear operators \oid{N}{}{reg}
as a "reflected" Grothendieck representation of the ideal $\oid{N}{}{}$. 
The Grothendieck representation states that $\oid{N}{}{} \stackrel{1}{=} 
\oid{W}{}{} \circ \oid{I}{}{}$ (cf. \cite{p1}, 24.6.2). 
  
\begin{theorem}
$(\oid{N}{}{d}, \idn{N}{}{d}) = (\oid{N}{}{reg}, \idn{N}{}{reg}) =
(\oid{I}{}{}, \idn{I}{}{}) \circ (\oid{K}{}{}, \| \cdot \|) =
(\oid{I}{}{}, \idn{I}{}{}) \circ (\oid{W}{}{}, \| \cdot \|)$.
\end{theorem}

{\noindent}{\sc PROOF}: Since $I$ is a (right-)accessible ideal it follows 
that $\oid{I}{}{} \circ \oid{G}{}{} \stackrel{1}{=} \oid{N}{}{}$, and lemma 
2.1 implies 
\begin{center}
$\oid{N}{}{reg} \stackrel{1}{=} (\oid{I}{}{} \circ \oid{G}{}{})^{reg}
\stackrel{1}{\subseteq} \oid{I}{}{} \circ \oid{G}{}{inj} 
\stackrel{1}{=} \oid{I}{}{} \circ \oid{K}{}{} \stackrel{1}{\subseteq}
\oid{I}{}{} \circ \oid{W}{}{}$.
\end{center} 
To prove the other inclusion observe that $\oid{I}{}{} \circ \oid{W}{}{}
\stackrel{1}{\subseteq} (\oid{W}{}{} \circ \oid{I}{}{})^{d}$ (since 
$\oid{W}{}{} \stackrel{1}{=} \oid{W}{}{d}$ and  
$\oid{I}{}{} \stackrel{1}{=} \oid{I}{}{d}$). Hence 
$\oid{I}{}{} \circ \oid{W}{}{} \stackrel{1}{\subseteq} 
\oid{N}{}{d} \stackrel{1}{=} \oid{N}{}{reg}$.
${}_{\displaystyle\Box}$\\

{\noindent}This representation leads to interesting consequences concerning
the ideal $\oid{K}{}{-1} \circ \oid{G}{}{}$, in particular to the conjugate
of this ideal. Note that $id_E \notin \oid{K}{}{-1} \circ \oid{G}{}{}$ for
each Banach space $E$ which has not the approximation property. To prepare 
this discussion we need the following
\begin{lemma}
$(\oid{G}{}{}, \| \cdot \|) \subseteq (\oid{I}{}{\Delta}, \idn{I}{}{\Delta})$.
\end{lemma}
{\noindent}{\sc PROOF}: Let $E, F$ be Banach spaces, $T\in \oid{G}{}{}(E, F)$ 
and  $L \in \oid{F}{}{}(F, E)$ be an arbitrary finite operator.
Then there exist $b_1, ...,b_n \in F'$ and $x_1, ...,x_n \in E$ such that
$Ly = \sum\limits_{i = 1}^n \langle y, b_i\rangle  x_i$ $(y \in F)$. Let 
$(T_m )_{m \in {\bf{N}}}$ be a sequence of finite operators such that 
$\lim\limits_{m\to\infty} \| T - T_m \| = 0$. Then for all $i \in \{ 1, ..., n \}$
$\lim\limits_{m\to\infty} \langle T_mx_i, b_i \rangle = \langle Tx_i, b_i \rangle$.
Thus 
\begin{center}
$|tr(TL)| = \lim\limits_{m\to\infty} |tr(T_mL)| \leq \| T \| \hspace{1mm} \idn{I}{}{}(L)$ 
\end{center}
which implies that $T \in \oid{I}{}{\Delta}(E, F)$ and $\idn{I}{}{\Delta}(T)
\leq \| T \|$.
${}_{\displaystyle\Box}$\\

{\noindent}{\bf{Remark:}} Let (\oid{A}{}{}, \idn{A}{}{}) an arbitrary 
$p$-Banach ideal $(0 < p \leq 1)$.
Using an analogous proof and the definition of \oid{A}{}{min} we obtain a 
generalization of the previous lemma:
\begin{center}
$(\oid{A}{}{min}, \idn{A}{}{min}) \subseteq 
(\oid{A}{}{{\displaystyle\ast}\Delta}, \idn{A}{}{{\displaystyle\ast}\Delta})$
${}_{\displaystyle\Box}$
\end{center}

\begin{prop}
$\oid{K}{}{-1} \circ \oid{G}{}{}$ is a totally accessible regular Banach
ideal which is not maximal. Moreover $\oid{K}{}{-1} \circ \oid{G}{}{} 
\stackrel{1}{\subseteq} \oid{N}{}{\Delta}$ and $\oid{I}{}{} \stackrel{1}{=} 
(\oid{K}{}{-1} \circ \oid{G}{}{})^{\Delta}$.
\end{prop}

{\noindent}{\sc PROOF}: Since $\oid{K}{}{}$ is an injective Banach ideal, 
theorem 2.1 implies that $\oid{N}{}{reg\Delta} 
\stackrel{1}{=} (\oid{I}{}{} \circ \oid{K}{}{})^{\Delta} \stackrel{1}{=}
\oid{K}{}{-1} \circ \oid{I}{}{\Delta}$. Hence $\oid{K}{}{-1} \circ 
\oid{G}{}{} \stackrel{1}{\subseteq} \oid{K}{}{-1} \circ \oid{I}{}{\Delta}
\stackrel{1}{\subseteq} \oid{N}{}{\Delta}$. Since $\oid{N}{}{\Delta}$ is
totally accessible it follows that (cf. \cite{oe3}, theorem 3.1.)
$\oid{I}{}{} \stackrel{1}{=} \oid{N}{}{\Delta\Delta} \stackrel{1}{\subseteq}
(\oid{K}{}{-1} \circ \oid{G}{}{})^{\Delta} \stackrel{1}{\subseteq} \oid{I}{}{}$. 
Suppose $\oid{K}{}{-1} \circ \oid{G}{}{}$ is a maximal Banach ideal. Then
$\oid{K}{}{-1} \circ \oid{G}{}{} \stackrel{1}{=} 
(\oid{K}{}{-1} \circ \oid{G}{}{})^{\displaystyle\ast\ast} \stackrel{1}{=}
(\oid{K}{}{-1} \circ \oid{G}{}{})^{\Delta\displaystyle\ast} \stackrel{1}{=}
\oid{I}{}{\displaystyle\ast} \stackrel{1}{=} \oid{L}{}{}$ which is a 
contradiction. The total accessibility of \oid{L}{}{} 
implies the total accessibility of $\oid{K}{}{-1} \circ \oid{G}{}{}$ and the 
regularity follows by  a straight forward calculation.  
${}_{\displaystyle\Box}$\\

{\noindent}We will finish this paper with another interesting application of
lemma 2.1.

\section{On operator ideals which factor through conjugates} 
In the following let $(\oid{A}{}{}, \idn{A}{}{})$ be an arbitrary maximal 
Banach ideal. Then the product ideal $\oid{A}{}{} \circ \oid{L}{\infty}{}$ 
is left-accessible (see \cite{oe3}, corollary 4.1). Using corollary 2.1 we will show 
that $\oid{A}{}{} \circ \oid{L}{\infty}{}$ is strongly related to 
(the right-accessible ideal) $\oid{A}{}{{\displaystyle\ast}\Delta}$ and 
to accessibility properties of the injective hull of 
$\oid{A}{}{\displaystyle\ast}$ (see \cite{oe3}, theorem 3.1).  
\begin{lemma}
Let $(\oid{A}{}{}, \idn{A}{}{})$ be a maximal Banach ideal. Then 
\begin{center}
$\oid{A}{}{} \circ\oid{L}{\infty}{} \stackrel{1}{\subseteq}
(\oid{A}{}{{\displaystyle\ast}\Delta} \circ \oid{L}{\infty}{})^{inj}
\stackrel{1}{\subseteq} (\oid{A}{}{{\displaystyle\ast}\Delta})^{inj}$.
\end{center} 
\end{lemma}

{\noindent}{\sc PROOF}: Let $E, F$ be Banach spaces, $\varepsilon > 0$ 
and $T \in \oid{A}{}{} \circ \oid{L}{\infty}{}(E, F)$. Then there exists 
a Banach space $G$, operators $R \in \oid{A}{}{}(G, F)$ and $S \in 
\oid{L}{\infty}{}(E, G)$ such that $j_FT = R''j_GS$ and $\idn{A}{}{}(R)\idn{L}
{\infty}{}(S) < (1 + \varepsilon)(\idn{A}{}{} \circ \idn{L}{\infty}{})(T)$.
Since $\oid{L}{\infty}{} \stackrel{1}{=} \backslash{\oid{L}{}{}}$ there are
a compact space $K$, operators $U \in \oid{L}{}{}(C(K), G'')$
and $V \in \oid{L}{}{}(E, C(K))$ such that $j_GS = UV$ and $\|U\|\|V\| <
(1 + \varepsilon) \idn{L}{\infty}{}(S)$. Hence $J_{F''}j_FT = (J_{F''}R''U)V$
- with canonical embedding $J_{F''} : F'' \stackrel{1}{\hookrightarrow} (F'')^{inj}$.
Since $R'' \in \oid{A}{}{}(G'', F'')$ and $\oid{A}{}{}(C(K), (F'')^{inj})
\stackrel{1}{=} \oid{A}{}{{\displaystyle\ast}\Delta}(C(K), (F'')^{inj})$
(both $C(K)$ and $(F'')^{inj}$ have the metric approximation property (cf.
\cite{oe1})) it follows that $J_{F''}j_FT \in 
(\oid{A}{}{{\displaystyle\ast}\Delta} \circ \oid{L}{\infty}{})(E, (F'')^{inj})$ 
and $(\idn{A}{}{{\displaystyle\ast}\Delta} \circ \idn{L}{\infty}{})(J_{F''}j_FT)
\leq \idn{A}{}{{\displaystyle\ast}\Delta}(J_{F''}R''U) \idn{L}{\infty}{}(V)
\leq \idn{A}{}{}(R) \hspace{1mm} \|U\| \|V\| < (1 + \varepsilon)\hspace{1mm} \idn{A}{}{}(R) 
\idn{L}{\infty}{}(S) < (1 + \varepsilon)^2 \hspace{1mm}
(\idn{A}{}{} \circ \idn{L}{\infty}{})(T)$.\\ 
Hence $T \in 
((\oid{A}{}{{\displaystyle\ast}\Delta} \circ \oid{L}{\infty}{})^{inj})^{reg}(E; F)$
and 
$((\idn{A}{}{{\displaystyle\ast}\Delta} \circ \idn{L}{\infty}{})^{inj})^{reg}(T)
\leq (\idn{A}{}{} \circ \idn{L}{\infty}{})(T)$.
Corollary 2.1 now yields the claim.
${}_{\displaystyle\Box}$\\

{\noindent}Note that \oid{A}{}{\displaystyle\ast inj} is left-accessible if
\oid{A}{}{} is {\it{injective}}. However lemma 3.1 and
(\cite{oe3}, proposition 4.1) imply a weaker condition:

\begin{prop}
Let \oid{A}{}{} be a maximal Banach ideal. If 
\oid{A}{}{{\displaystyle\ast}\Delta} is injective then 
\oid{A}{}{\displaystyle\ast inj} is (totally) accessible.
\end{prop}

{\noindent}{\sc PROOF}: Since \oid{A}{}{{\displaystyle\ast}\Delta} is
injective, lemma 3.1 implies that 
$\oid{A}{}{} \circ \oid{L}{\infty}{} \stackrel{1}{\subseteq} 
\oid{A}{}{{\displaystyle\ast}\Delta}$. Hence 
$\oid{A}{}{\displaystyle\ast} \circ \oid{A}{}{} \circ \oid{L}{\infty}{} 
\stackrel{1}{\subseteq} \oid{I}{}{}$ and it follows that $\oid{A}{}{\displaystyle\ast} 
\circ \oid{A}{}{} \stackrel{1}{\subseteq} \oid{P}{1}{}$. 
(\cite{oe3}, proposition 4.1) finishes the proof.
${}_{\displaystyle\Box}$\\

\end{document}